\documentclass[11pt]{article}   
\usepackage{amssymb,amscd,latexsym}   
\usepackage{amsmath}
\usepackage{amsthm}
\newtheorem{Theorem}{Theorem}[section]
\newtheorem{Lemma}[Theorem]{Lemma}

\newtheorem{Proposition}[Theorem]{Proposition}
\newtheorem{Remark}[Theorem]{Remark}
\newtheorem{Example}[Theorem]{Example}

\newtheorem{Definition}[Theorem]{Definition}
\newtheorem{Question}[Theorem]{Question}
\newcommand{\rar}{\rightarrow}
\newcommand{\lar}{\longrightarrow}
\newcommand{\llar}{-\kern-5pt-\kern-5pt\longrightarrow}

\newcommand{\injects}{\hookrightarrow}

\renewcommand{\phi}{\varphi}
\newcommand{\demo}{{\sc Proof. }}
\renewcommand{\proof}{\demo}

\newcommand{\xx}{{\bf x}}
\newcommand{\yy}{{\bf y}}

\newcommand{\ff}{{\bf f}}
\renewcommand{\gg}{{\bf g}}

\newcommand{\hht}{{\rm ht}\,}

\newcommand{\restr}{{\kern-1pt\restriction\kern-1pt}}

\newcommand{\pp}{{\mathbb P}}

\begin{document}

\begin{center}
{\Large{\bf\sc New constructions of Cremona maps}}
\footnotetext{2010 AMS {\it Mathematics Subject
Classification}: 13D02, 13H10, 14E05, 14E07, 14M05, 14M25.}

\vspace{0.3in}

{\large\sc Barbara Costa}      \quad\quad
 {\large\sc Aron  Simis}\footnote{Partially
supported by a CNPq grant.}

\end{center}

\begin{abstract}

One defines two ways of constructing rational maps
derived from other rational maps, in a characteristic-free context.
The first introduces the Newton complementary dual of a rational map.
One main result is that this dual preserves birationality and gives an
involutional map of the Cremona group to itself that restricts to the monomial
Cremona subgroup and preserves de Jonqui\`eres maps. In the monomial restriction this duality
commutes with taking inverse in the group, but is a not a group homomorphism.
The second construction is an iterative process to obtain rational maps in increasing
dimension. Starting with birational maps, it leads to rational maps whose topological degree
is under control. Making use of monoids, the resulting construct is in fact birational if the
original map is so.
A variation of this idea is considered in order to preserve properties of the base ideal,
such as Cohen--Macaulayness. Combining the two methods, one is able to produce explicit
infinite families of Cohen--Macaulay Cremona maps with prescribed
dimension, codimension and degree.

\end{abstract}


\section*{Introduction}
\label{intro}

The purpose of this paper is to introduce a couple of constructions of rational maps out of given ones.
The core of the results consists in that the constructions preserve birationality and often a few other properties.
The source and the target of any rational map considered here are
projective spaces over a field $k$. Such a map will typically be denoted $\mathfrak{F}:\pp^n\dasharrow \pp^m$,
where the field $k$ is self-understood.
For the geometric purpose $k$ is usually taken to be algebraically closed. However, we stress that the constructions
and the main theorems are characteristic-free.

The first of these constructions is a large extension of an operation introduced in
\cite{birational-linear} for rational maps defined by monomials. It associates to $\mathfrak{F}:\pp^n\dasharrow \pp^m$ another rational
map with same source and target which will be birational (onto the image) if $\mathfrak{F}$ is so.
In particular, if $\mathfrak{F}$ is a Cremona map of $\pp^n$, so is its associated  construct.
Since the set of Cremona maps of $\pp^n$ form a group under composition, it is obvious that one can produce
new Cremona maps galore  out of given ones by means of group theoretic operations, such as generation and conjugation.
However, the present construction is of a different nature and resembles to a ``duality''  on the Cremona
group -- in that applying twice successively yields back the original map, -- and to a group homomorphism, when
restricted to the monomial Cremona subgroup -- in that it preserves inverses (but not the group operation).

The second construction is an inductive procedure associating to $\mathfrak{F}:\pp^n\dasharrow \pp^m$ a rational
map  $\mathfrak{G}:\pp^{n+1}\dasharrow \pp^{m+1}$.
A germ of this idea has been given in \cite{Barbara}, but it is possible that examples of
this construction have been considered in the classical references.
Nonetheless, our interest is as to how the properties of the given rational map transfer over
to the construct.
The procedure preserves birationality and, in addition, properties of the base ideal.
A main purpose of this construction here is to explicitly produce, for arbitrarily given integers $n\geq 2$,
$d\geq 1$ and $r$ in the integer interval $[2,n]$,  infinitely many
Cremona maps of $\pp^n$ of degree
$d\geq 1$ and  codimension $r$; in addition, if $d\geq n+1-r$, there are infinitely many such maps whose base ideal
is Cohen--Macaulay.
To our knowledge, such families do not seem to have been explicitly stated in the literature.

\medskip

Our central results are stated as Theorem~\ref{dual_forms_is_birational},
Theorem~\ref{CM_birational_cod2}, and Theorem~\ref{arbitrary_cod}.
They are better understood within the context of the respective section.
We accordingly give a brief description of the corresponding sections.

Section~\ref{sec:1} starts with a swift recap of birational maps, after which we
introduce the concept of the Newton matrix $N(\ff)$ of a set $\ff=\{f_0,\ldots,f_m\}$ of forms of the same degree
in the standard graded polynomial ring $R=k[x_0,\ldots,x_n]$ over a field $k$.
This matrix is briefly described as the concatenation of the log-matrices of all
terms appearing in these forms.
One should be careful not to abuse too much handling the ideals generated by the forms
or its various parts, as it is rather the subalgebra generated by those that will
play a distinctive role.
The subsequent main construct is what one calls the Newton complementary dual matrix $\widehat{N(\ff)}$
of the Newton matrix $N(\ff)$. This in turn gives rise, in a sort of ``coefficient-frame'' operation,
to a new set of forms $\widehat{\ff}$ in the same ring $R$, called similarly
the Newton complementary dual set of $\ff$.
From the definitions follows that $\widehat{\widehat{\ff}}=\ff$.
The first relevant result is the statement that if $\ff$ are monomials, then $\ff$ defines a birational map
(onto the image) if and only if $\widehat{\ff}$ does too.
The proof is an elegant field-theoretic argument based on the characteristic property of the construction
in the monomial case (see Remark~\ref{achtung} (3)).
In particular, it restricts to a well-defined ``duality'' of the monomial subgroup of the entire Cremona group
of dimension $n$.
Another meaningful property of this duality in the case of monomials is that it commutes with taking
inverses.
However, this property no longer subsists for arbitrary birational maps, as is shown by examples.

At any rate, the main point of this section is proving that birationality is preserved for arbitrary forms
and, in addition, there is a $k$-isomorphism of the graded $k$-subalgebras $k[\ff]$ and $k[\widehat{\ff}]$.
As a non-trivial illustration of the theory, we show that Newton complementary duality preserves plane de Jonqui\`eres
maps. However, the degree of the resulting map is given by a rather involved expression in terms of the defining
forms of the original map -- perhaps not too surprising, as the degree of a plane Cremona map is dependent upon
the multiplicities of the infinitely near base points of the map.
Unfortunately, Newton complementary duality does not preserve the composition of Cremona maps.
Therefore, it is less natural to ask what group theoretic properties are preserved by this construction.

\smallskip

In the second section we deal with a procedure of associating to a  rational map $\mathfrak{F}:\pp^n\dasharrow \pp^m$ a rational
map  $\mathfrak{G}:\pp^{n+1}\dasharrow \pp^{m+1}$, by viewing $\pp^n$ as a cone and letting the last coordinate of $\mathfrak{G}$
be a monoid with respect to the cone variable.
This will have the effect that the image of $\mathfrak{G}$  will be a cone over the image of $\mathfrak{F}$.
Such rigid features of the construction will entail the preservation of birationality and typical properties of the
base ideal.

In a more explicit way, let $\ff\subset R=k[x_0,\ldots, x_{n-1}]$ be forms of degree $d$ defining $\mathfrak{F}$.
Pick $g=g_{d-1}x_n+g_d\in R[x_n]=k[x_0,\ldots, x_n]$, where $g_{d-1}\neq 0$ and $g_d$ are forms in $R$ of degrees $d-1$ and $d$,
respectively. Then the forms $\ff,g$ define $\mathfrak{G}$.
We show that $\mathfrak{G}$ is birational along with $\mathfrak{F}$.
Moreover, if $g_d\in I$ and $g_{d-1}$ is a non-zero-divisor on $R/I$ then the Cohen--Macaulay property transfers as well.
This is obtained as a particular case of a more general principle when one takes $g$ to be a form of higher degree in
the new variable $x_n$, by which the topological degree $\mathfrak{G}$ may turn out to be non-trivial depending on that degree, even
when $\mathfrak{F}$ is birational.
In the sequel, we give a slight variation of this iterative procedure in order to preserve the condition
that the base ideal is Cohen--Macaulay of codimension $2$.
In this form, the construction resembles the de Jonqui\`eres parametrizations of \cite{stellar}.
However, the resemblance is superficial since there the purpose is implicitization, while here the main use is in the Cremona case.
Combining the two methods developed in this section, we are able to prove the existence of infinite families
of Cohen--Macaulay Cremona tranformations with prescribed dimension, codimension and degree (bounded below in a natural
way).

At the end of section we show how to produce codimension $2$ Cohen--Macaulay  Cremona maps that are very close
to the determinantal polar maps introduced in \cite{CRS}.
These examples may bear a meaning for the theory of determinantal polar maps, a rather narrow habitat.

\section{The Newton complementary dual of a rational map}
\label{sec:1}

\subsection{Recap of birational maps}
\label{sec:1.1}

Our reference for the basics in this part is \cite{bir2003}, which contains enough of the introductory
material in the form we use here (see also \cite{AHA} for a more general overview).

Let $k$ denote an arbitrary infinite field which will be assumed to be algebraically closed for
the geometric purpose.
A rational map $\mathfrak{F}:\pp^n\dasharrow \pp^m$ is defined by $m+1$ forms $\ff=\{f_0,\ldots, f_m\}
\subset R:=k[\xx]=k[x_0,\ldots,x_n]$ of the same degree $d\geq 1$, not all null.
We often write $\mathfrak{F}=(f_0:\cdots :f_m)$ to underscore the projective setup.
Any rational map can without lost of generality be brought to satisfy  the condition
that $\gcd\{f_0,\cdots ,f_m\}=1$ (i.e., the linear system spanned by $\ff$
 {\em has no fixed part}).
One also assumes that every variable $x_i$ appears effectively in some $f_j$.
If the rational map satisfies both conditions we will say that is satisfies the {\em canonical restrictions}.
These conditions are usually hand-waved when considering Cremona maps of dimension $n$.
The reason is that one can rather take equivalence classes of sets $\ff$ as above,
of the same cardinality, by multiplication with an arbitrary form; clearly, every such class has a unique representative
satisfying $\gcd\{f_0,\cdots ,f_m\}=1$. Composition in the Cremona
group is well-defined for these classes, which allows for the slight get-away of working with these unique
representatives (see \cite{AHA} for a detailed account of such preliminaries, but of course this is well-known in
the classical subject).

The common degree $d$ of the $f_j$ is the {\em degree} of $\mathfrak{F}$, not to be confused with its
topological (or field) degree.
The ideal $I_{\mathfrak{F}}:=(f_0,\ldots,f_m)\subset R$ is called the {\em base ideal} of $\mathfrak{F}$.
The {\em image} of $\mathfrak{F}$ is the projective subvariety $W\subset \pp^m$ whose homogeneous
coordinate ring is the $k$-subalgebra $k[\ff]\subset R$ after degree renormalization.
Write $S:=k[\ff]\simeq k[\yy]/I(W)$, where $I(W)\subset k[\yy]=k[y_0,\ldots,y_m]$ is the homogeneous defining ideal
of the image in the embedding $W\subset \pp^m$.

We say that $\mathfrak{F}$ is {\em birational onto the image} if there is a rational map
backwards $\pp^m\dasharrow \pp^n$ such that the residue classes $\ff'=\{f'_0,\ldots, f'_n\}
\subset S$ of a set of defining coordinates do not simultaneously vanish and satisfy the
relations
\begin{equation}\label{birational_rule}\nonumber
(\ff'_0(\ff):\cdots :\ff'_n(\ff))=(x_0:\cdots :x_n), \;
(\ff_0(\ff'):\cdots :\ff_m(\ff'))\equiv (y_0:\cdots :y_m)\pmod {I(W)}
\end{equation}
Let $K$ denote the field of fractions of $S=k[\ff]$.
The coordinates $\{f'_0,\cdots ,f'_n\}$ defining the ``inverse'' map are not uniquely defined;
any other set $\{f''_0,\cdots ,f''_n\}$ related to the first through requiring that it define the same
element of the projective space $\pp^n_{K}=\pp^n_k\otimes_k {{\rm Spec}(K)}$ will do as well -- both tuples are called
{\em representatives} of the rational map (see \cite{bir2003} for details).
One can see that, if $k$ is algebraically closed, these relations
translate into the geometric definition in terms of invertibility of the map on a dense Zariski open
set.

An alternative definition of a birational map is in terms of field extensions.
Letting the common degree of the $f$'s be $d$, one has an inclusion of rings
$S=k[\ff]\subset R^{(d)}:=k[\xx_d]$. Then $\mathfrak{F}$ is birational onto the image if and only
if this inclusion triggers an equality $K=k(\xx_d)$ at the level of the respective fields of fractions.
The inverse map is given by extracting one representative from the inverse
field isomorphism.

\subsection{The Newton complementary dual}
\label{sec:1.2}

The following notion extends to arbitrary forms of the same degree in $R=k[x_0,\ldots,x_n]$
the concept of a certain dual construction of an integer matrix (see \cite{Barbara}, \cite{escobar},
\cite{birational-linear}).

\begin{Definition}\label{def_of_dual}\rm
Let $f$ be $d$-form in $R$, for some $d\geq 1$.
Denote by $N(f)$ (as a reminder of the Newton polygon) the so-called log-matrix of the set of nonzero terms
of $f$, where we assume that the terms are lexicographically ordered.
Thus, $N(f)$ is the matrix whose columns are the
exponents vectors of the nonzero terms of $f$ lexicographically ordered.
One may call $N(f)$ the {\sc Newton log matrix} (or simply the {\sc Newton matrix}) of $f$.
\end{Definition}
Given an ordered set $\mathbf{f}:=\{f_0,\ldots,f_m\}$ of such forms of the same degree $d\geq 1$,
let $N(\mathbf{f})$ denote the concatenation of the Newton matrices $N(f_0),\ldots,N(f_m)$;
accordingly, we call $N(\mathbf{f})$ the {\sc Newton matrix} of the set $\mathbf{f}$.
Note that $N(\mathbf{f})$ is an integer stochastic matrix.

The row vector $\mathbf{c}_f$ whose entries are the nonzero coefficients of a form $f$, ordered by the lexicographic monomial order, is
the {\em coefficient frame} of $f$.
We write symbolically
\begin{equation}\label{symb}
f=\langle\mathbf{c}_f,\xx^{N(f)}\rangle
\end{equation}
as the inner product of the coefficient frame by the set of the corresponding monomials.

The {\sc Newton complementary dual matrix} (or simply the {\sc complementary dual matrix}) of
$N(\ff)=(a_{i,\ell})$ is the matrix
$$\widehat{N(\ff)}=(\alpha_i - a_{i,\ell}),$$
where, for every $0\leq i\leq n$, $\alpha_i = \max_{\ell} \{a_{i,\ell}\}$,  with $\ell$ indexing the set of all nonzero terms of all
forms in the set $\mathbf{f}$.

In other words, denoting $\underline{\alpha}:=(\alpha_1,\ldots,\alpha_n)^t$, one has
$$
\widehat{N(\ff)}=\left[\,\underline{\alpha}\,| \cdots |\,\underline{\alpha}\,\right]_{n \times (r_0+\ldots + r_m)} - N(\ff),
$$
where $r_j$ denotes the number of nonzero terms of $f_j$, $j=0,\ldots, m$.

For every $j=0,\ldots,m$, let $\widehat{N(\ff)}_j$ denote the submatrix of $\widehat{N(\ff)}$
whose columns come from $f_{j}$.
Finally consider the set of forms
\begin{equation}\label{dual_forms}
\widehat{\ff}:=\{\widehat{f_0}:=\langle\mathbf{c}_0,\xx^{\widehat{N(\ff)}_0}\rangle\,,\ldots, \,\widehat{f_m}:=
\langle\mathbf{c}_m,\xx^{\widehat{N(\ff)}_m}\rangle\}
\end{equation}
where $\mathbf{c}_j=\mathbf{c}_{f_j}$ stands for the coefficient frame of $f_j$.

\begin{Remark}\label{achtung}\rm
(1) It is extremely important to note that the submatrix $\widehat{N(\ff)}_j$ is in general different from
the Newton complementary dual matrix $\widehat{N(f_j)}$ of the set consisting solely of $f_j$.
In this regard the above notation $\widehat{f_j}$ may lead to confusion, but it will always mean the $j$th element of
the set $\widehat{\ff}$ while applying
the Newton procedure to the whole set $\ff$ and not step by step to its individual elements.
The process is always controlled by the ``global'' vector $\underline{\alpha}$ above which is clearly different
for the whole of $\ff$ as for one of its subsets.

(2) By the same token, the notion is not ``recurrent'' in the sense that we cannot in general  inductively reconstruct
$\widehat{\ff}$ from $\widehat{\gg}$, where $\gg=\{f_0,\ldots,f_{m-1}\}$.
Actually, without further restrictions the procedure can lead to very uninteresting situations; for example, if $\ff=\{x_0\}$
then $\widehat{\ff}=\{1\}$ !

(3) If $M$ is a monomial in the list of nonzero terms of $\ff$ then $M\cdot\widehat{M}=\underline{x}^{\alpha}$,
with  $\underline{x}^{\alpha}=x_0^{\alpha_0}\cdots x_n^{\alpha_n}$, where the $\alpha$'s are as above.

(4) The common degree of the forms in $\widehat{\ff}$ is the {\em dual degree} $\hat{d}:=\alpha_0+\cdots +\alpha_n - d$.
Note the immediate bounds $n+1-d\leq \hat{d}\leq nd$, where the first inequality follows from the non-degeneracy
of the forms in $\mathbf{f}$ (and  tells us literally nothing when $d>n+1$).
Anyway, both are clearly easily attainable in the family of arbitrary sets $\ff$.
However, for special relevant subfamilies these bounds can be tightened.
Thus, for instance, if $\ff$ defines a Cremona map of $\pp^n$ of degree $d\geq 2$ then a slightly better upper bound is $nd-1$.
It is challenging to lower this upper bound for families of Cremona maps of $\pp^n$ enjoying special properties.
\end{Remark}

\begin{Definition}\rm
We will call $\widehat{\ff}$ the {\sc Newton complementary dual} set of $\mathbf{f}$.
\end{Definition}

The terms ``Newton'' and ``complementary'' seem to be appropriate.
The following easy result justifies the term ``dual''.

\begin{Lemma}\label{minimax}
If $\,\ff$ satisfies the canonical restrictions  then so does $\widehat{\ff}$ and one has $\widehat{\widehat{\ff}}=\ff$.
\end{Lemma}
\proof
Since $\ff$ has no fixed part, in particular no variable belongs to the support of every term of $\ff$.
Therefore, for any given $i=0,\ldots, n$, some $a_{il}=0$.
On the other hand, $\ff$ has also the property that for every $i=0,\ldots, n$, some $a_{ik}\neq 0$.
Collecting the two gives that for any $i$, $\alpha_i\geq 1$.
We may assume that this is the maximum.
It follows that $\alpha_i-a_{ik}=0$. We thus conclude that $\widehat{\ff}$ satisfies the canonical restrictions.

Next, write $\widehat{N(\ff)}=(b_{ij})$.
By definition, $b_{ij}=\alpha_i-a_{ij}$, for every $j=0,\ldots,n$, where
$\alpha_i = \max_j \{a_{ij}\}$.
Let $\beta_i=\max_{j}\{b_{ij}\}=\max_{j}\{\alpha_i-a_{ij}\}$.
Then $\beta_i=\alpha_i-\min_j\{a_{ij}\}=\alpha_i$, as $\min_j\{a_{ij}\}=0$. Thence
$$(\widehat{\widehat{N(\ff)}})_{ij}=\beta_i-b_{ij}=\alpha_i-(\alpha_i-a_{ij})=a_{ij},$$
the last equality following from the hypothesis that $\ff$ has no fixed part.
Therefore, $\widehat{\widehat{N(\ff)}})=N(\ff)$, thus showing the contention.
\qed

\smallskip

We will call the vector $\underline{\alpha}:=(\alpha_1,\ldots,\alpha_n)^t$ the {\sc directrix vector}
of $\ff$.
 Note that the directrix vectors of $\ff$ and of $\widehat{\ff}$ coincide, as has been seen in the proof
 of the lemma.

\begin{Example}\rm If $\ff=\{x_0,\ldots, x_n\}$ then its directrix vector is $(1,\ldots,1)^t$ and
$$\widehat{\ff}=\{x_1\cdots x_n,\, \ldots, \, x_0\cdots \hat{x_i}\cdots x_n,
\,\ldots, \,x_0\cdots x_{n-1}\}.$$
In terms of the rational maps defined by the respective forms, this says that the Newton complementary dual of the identity
map of $\pp^n$ is the Magnus reciprocal involution.
\end{Example}

Note the relations $f_j\widehat{f_j}=x_0\cdots x_n$ for every $j=0,\ldots, n$.
The example, along with this property, is  propaedeutic for the result in the next section.

\subsection{The main result}
\label{sec:1.3}

For convenience, we often say that a set $\ff$ of forms of the same degree satisfying the canonical restrictions
is birational if the associated rational map is birational onto its image.
Note that this notion is independent of the ordering of the forms.

We first deal with the case where the forms in $\ff$ are monomials.

\begin{Proposition}\label{dual_mons_is_birational}
Let $\ff=\{f_0,\ldots,f_m\}\subset R=k[x_0,\ldots,x_n] (n\geq 1)$ be monomials of the same degree satisfying
the canonical restrictions.
If $\ff$ defines a birational map onto its image then so does its
Newton complementary dual set $\widehat{\ff}$.
\end{Proposition}
\proof Let $d$ be the common degree of the monomials in $\ff$.

Note that $k(\underline{x}_d)=k(x_0^d,x_1/x_0,\ldots,x_n/x_0)$, hence the birationality of $\ff$ is tantamount to the equality
$$k(f_0,\ldots,f_m)=k\left(x_0^d,\frac{x_1}{x_0},\ldots,\frac{x_n}{x_0}\right).
$$
Therefore, for every $i=1,\ldots,n$ one has an expression
$$\frac{x_i}{x_0}=\frac{g_i(f_0,\ldots,f_m)}{g'_i(f_0,\ldots,f_m)},$$
where $g_i, g'_i$ are forms in the polynomial ring $k[y_0,\ldots,y_m]$ of the same degree, say, $s$.

By the same token, it suffices to show the equality
$$k(\widehat{f_0},\ldots,\widehat{f_m})
=k\left(x_0^{|\underline{\alpha}|-d},\frac{x_1}{x_0},\ldots,\frac{x_n}{x_0}\right)\;(=k(\underline{x}_{|\underline{\alpha}|-d})),$$
where $\widehat{\ff}=\{\widehat{f_0},\ldots,\widehat{f_m}\}$ and $\underline{\alpha}$ is the directrix vector and
$|\underline{\alpha}|=\alpha_0+\cdots +\alpha_m$.
Once more we insist on the meaning of $\widehat{f_j}$ as explained at the beginning of the section.

Now, since the forms in $\ff$ are monomials, one has $\widehat{f_j}=\underline{x}^{\alpha}/f_j$ for every $j=0,\ldots, m$,
where  $\underline{x}^{\alpha}=x_0^{\alpha_0}\cdots x_n^{\alpha_n}$.
It follows that, for every $i\in \{1,\ldots,n\}$, one has
$$\frac{x_i}{x_0}=\frac{g_i(f_0,\ldots,f_m)}{g'_i(f_0,\ldots,f_m)}\frac{(\underline{x}^{\alpha})^s}{(\underline{x}^{\alpha})^s}
= \frac{g_i(f_0,\ldots,f_m)}{(\underline{x}^{\alpha})^s}\frac{1}{\frac{g'_i(f_0,\ldots,f_m)}{(\underline{x}^{\alpha})^s}}=$$
$$\frac{g_i(\frac{1}{\widehat{f_0}},\ldots,\frac{1}{\widehat{f_m}})}{g'_i(\frac{1}{\widehat{f_0}},\ldots,\frac{1}{\widehat{f_m}})}
=\frac{h_i(\widehat{f_0}, \ldots, \widehat{f_m})}{h'_i(\widehat{f_0}, \ldots, \widehat{f_m})},$$
for suitable forms $h_i,h'_i\in k[y_0,\ldots,y_m]$ of the same degree.
It remains to deal with the power $x_0^{|\underline{\alpha}|-d}$, for which we claim likewise that
it belongs to $k(\widehat{f_0},\ldots,\widehat{f_m})$.
Since $\ff$ is a set of monomials satisfying the canonical restrictions, there is an index $j$ such that
$x_0$ does not divide $f_j$. Without loss of generality, say, $j=0$.
Write accordingly $f_0= x_{i_1} \cdots x_{i_d}$ with $0< i_1 \leq i_2 \leq \ldots \leq i_d$.
Using the above expressions, one has
$$\frac{f_0}{x_0^d}=\frac{x_{i_1}}{x_0} \cdots \frac{x_{i_d}}{x_0}=
\frac{h_{i_1}(\widehat{f_0}, \ldots, \widehat{f_m})}{h'_{i_1}(\widehat{f_0}, \ldots, \widehat{f_m})} \cdots
\frac{h_{i_d}(\widehat{f_0}, \ldots, \widehat{f_m})}{h'_{i_d}(\widehat{f_0}, \ldots, \widehat{f_m})}.$$
On the other hand,
\begin{eqnarray*}\nonumber
\frac{f_0}{x_0^d}&=&\frac{f_0}{x_0^d}\frac{x_0^{|\underline{\alpha}|}}{x_0^{|\underline{\alpha}|}}
=x_0^{|\underline{\alpha}|-d} \frac{f_0}{x_0^{|\underline{\alpha}|}}
=x_0^{|\underline{\alpha}|-d} \frac{f_0}{\underline{x}^{\alpha}} \frac{\underline{x}^{\alpha}}{x_0^{\alpha_0} \cdots x_0^{\alpha_n}}
=x_0^{|\underline{\alpha}|-d} \frac{1}{\widehat{f_0}} \,\frac{x_1^{\alpha_1} \ldots x_n^{\alpha_n}}{x_0^{\alpha_1} \ldots x_0^{\alpha_n}}\\
&=& x_0^{|\underline{\alpha}|-d} \frac{1}{\widehat{f_0}} \left(\frac{h_{i_1}(\widehat{f_0}, \ldots, \widehat{f_m})}
{h'_{i_1}(\widehat{f_0}, \ldots, \widehat{f_m})}\right)^{\alpha_1} \cdots \;\left(\frac{h_{i_d}(\widehat{f_0}, \ldots, \widehat{f_m})}
{h'_{i_d}(\widehat{f_0}, \ldots, \widehat{f_m})}\right)^{\alpha_n}.
\end{eqnarray*}
Equating the above two expressions yields the required result.
\qed

\medskip

Besides preserving birationality in the case of monomials, the dual process has an additional crucial purely algebraic property:

\begin{Lemma}\label{same_algebra}
Let $\ff=\{f_0,\ldots,f_m\}\subset R=k[x_0,\ldots,x_n] (n\geq 1)$ be monomials of the same degree  satisfying
the canonical restrictions, and let $\widehat{\ff}$ denote its Newton complementary dual set.
Then the identity map of $k[y_0,\ldots,y_m]$ induces a $k$-algebra isomorphism $k[\ff]\simeq k[\widehat{\ff}]$.
\end{Lemma}
\proof
Recall the basic relationship between the two in the monomial case: $f_j=\xx^{\alpha}/\widehat{f_j}$, for $j=0,\ldots,m$,
where $\alpha$ gives the coordinates of the directrix vector $\underline{\alpha}$.

It suffices to show that any polynomial relation of $\ff$ is one of $\widehat{\ff}$, and conversely too.
Since both algebras are toric and homogeneous, it suffices to consider homogeneous binomial relations.
Thus, let $\yy^{\beta}-\yy^{\gamma}\in k[y_0,\ldots,y_m]$, with $|\beta|=|\gamma|$.
This binomial is a relation of $\ff$ if and only if any one of the following equivalent conditions hold:
\begin{eqnarray*}\nonumber
\prod_j f_j^{\gamma_j}\kern-5pt&=&\kern-3pt\prod_j f_j^{\beta_j} \Leftrightarrow
\frac{(\underline{x}^{\alpha})^{|\gamma|}}{\prod_j f_j^{\gamma_j}} =
\frac{(\underline{x}^{\alpha})^{|\beta|}}{\prod_j f_j^{\beta_j}}
\Leftrightarrow \frac{\prod_j (\underline{x}^{\alpha})^{\gamma_j}}{\prod_j f_j^{\gamma_j}} =
\frac{\prod_j (\underline{x}^{\alpha})^{\beta_j}}{\prod_j f_j^{\beta_j}}
\Leftrightarrow \prod_j (\frac{\underline{x}^{\alpha}}{f_j})^{\gamma_j} =
\prod_j (\frac{\underline{x}^{\alpha}}{f_j})^{\beta_j}\\
&\Leftrightarrow& \prod_j (\widehat{f_j})^{\gamma_j} = \prod_j (\widehat{f_j})^{\beta_j}.
\end{eqnarray*}
Since the last equality above means that the binomial $\yy^{\beta}-\yy^{\gamma}$ is a relation of
$\widehat{\ff}$, we are through.
\qed

\begin{Remark}\rm
Note that the basic relations $f_j=\xx^{\alpha}/\widehat{f_j}$, for $j=0,\ldots,m$, imply an isomorphism
of the algebras $k[\ff]$ and $k[1/\widehat{f_0},\ldots, 1/\widehat{f_m}]$.
However, for arbitrary forms $g_0,\ldots,g_m$ of the same degree, there  may not be an isomorphism between $k[g_0,\ldots,g_m]$
and $k[1/g_0,\ldots,1/g_m]$ induced by the identity map of $ k[y_0,\ldots,y_m]$.
An example is provided by noting that the Grassmann--Pl\"ucker quadratic relation of the $2$-minors of a $2\times 4$
matrix is not a relation of the reciprocals of these minors.
\end{Remark}

One is now ready for the main result of this part.

\begin{Theorem}\label{dual_forms_is_birational}
Let $\ff=\{f_0,\ldots,f_m\}\subset R=k[x_0,\ldots,x_n] (n\geq 1)$ be arbitrary forms of the same degree satisfying
the canonical restrictions and let $\widehat{\ff}$ denote its Newton complementary dual set.
One has:
\begin{enumerate}
\item[{\rm (a)}] For suitable indeterminates $y_0,\ldots,y_m$ over $k$, the identity map of $k[y_0,\ldots,y_m]$
induces an isomorphism $k[\ff]\simeq k[\widehat{\ff}]$.
\item[{\rm (b)}] If $\ff$ defines a birational map onto its image then so does  $\widehat{\ff}$.
\end{enumerate}
\end{Theorem}
\proof
Keeping the notation of (\ref{symb}), write $f_j=\langle\mathbf{c}_{f_j},\xx^{N(f_j)}\rangle$, for $j=0,\ldots, m$,
where $\mathbf{c}_{f_j}=(c_{j,1},...,c_{j,r_j})$, so that $f_j=\langle\mathbf{c}_{f_j},\xx^{N(f_j)}\rangle=
\sum_{\ell=1}^{r_j} c_{j,\ell}\,M_{j,\ell}$, for suitable monomials $M_{j,\ell}$.
In this notation, one has $\widehat{f_j}=\sum_{\ell=1}^{r_j} c_{j,\ell}\,\widehat{M_{j,\ell}}$.
Let $\{z_{j,1},\ldots, z_{j,r_j}\}$, with $j=0,\ldots, m$, denote $r_0+\cdots +r_m$ mutually independent indeterminates,

By Lemma~\ref{same_algebra}, the identity map of $k[z_{j,1},\ldots, z_{j,r_j}\,| 0\leq j\leq m]$ induces an isomorphism
of $k[\mathbf{M}]:=k[M_{j,1},\ldots, M_{j,r_j}\,| 0\leq j\leq m]$ onto $k[\widehat{\mathbf{M}}]:=
k[\widehat{M_{j,1}},\ldots, \widehat{M_{j,r_j}}\,|\, 0\leq j\leq m]$.
Setting $y_j=\sum_{\ell=1}^{r_j} c_{j,\ell}\,z_{j,\ell}$, for  $0\leq j\leq m$,  yields algebraically independent
elements $y_0,\ldots,y_m$ over $k$.
Clearly, the restriction of the above map gives that the identity map of the polynomial ring $k[y_0,\ldots,y_m]$
induces an isomorphism of $k[\ff]$ onto $k[\widehat{\ff}]$.
This shows part (a).

Part (b) goes as follows: by part (a), passing to the respective fields of fractions on both sides, yields
a diagram of isomorphisms and inclusions:
$$\begin{array}{ccc}
k(\xx_d) & & k(\xx_{|\underline{\alpha}|-d})\\
\cup & & \cup\\
k(\mathbf{M}) &\simeq & k(\widehat{\mathbf{M}})\\
\cup &  &\cup\\
k(\ff) & \simeq & k(\widehat{\ff})
\end{array}
$$
where $\alpha$ was explained earlier.
The birationality assumption on $\ff$ implies that  $k(\ff)=k(\xx_d)$. In particular, $k(\ff)=k(\mathbf{M})$.
Since all horizontal maps are induced by the identity map on the polynomial ring $k[z_{j,1},\ldots, z_{j,r_j}\,| 0\leq j\leq m]$,
it follows that $k(\widehat{\ff})=k(\widehat{\mathbf{M}})$.
Finally, by the monomial case in Lemma~\ref{dual_mons_is_birational}, one has $k(\widehat{\mathbf{M}})=k(\xx_{|\underline{\alpha}|-d})$.
Therefore, $k(\widehat{\ff})=k(\xx_{|\underline{\alpha}|-d})$, as was required to show.
\qed

\subsection{Notable Cremona maps}
\label{sec:1.4}

\subsubsection{Monomial Cremona maps}

The operations of taking the
Newton dual and the inverse map commute in the case of monomial Cremona maps, according to  \cite[Theorem 3.1.3]{Barbara}.
For convenience we restate this result along with its proof, which depends on \cite[Theorem 2.2]{CremonaMexico}
(see also \cite[Theorem 2.2]{CostaSimis} for a restatement stressing the canonical restrictions).

For light reading, if $\ff$ defines a Cremona map, we let $\ff^{-1}$ denote its inverse map.
Recall that if $\ff$ is a set of monomials of the same degree, its Newton matrix is just the log-matrix $L(\ff)$ of these
monomials.
In the case of a monomial Cremona map of $\pp^n$, there is a meaningful vector
of $\mathbb{N}^{n+1}$ called the {\em inversion vector} (see \cite[Section 2]{CostaSimis}).
Along with the directrix vector $\underline{\alpha}^t\in \mathbb{N}^{n+1}$ it constitutes a key for many structural
results in monomial Cremona theory.

\begin{Proposition}\label{operations_commute}
Let $\ff\subset R=k[x_0,\ldots,x_n]$ be a set of $n+1$ monomials of the same degree satisfying
the canonical restrictions and defining a Cremona map of $\pp^n$.
Then $(\widehat{\ff})^{-1}$ and $\widehat{{\ff}^{-1}}$ coincide as maps.
More exactly,
\begin{equation}\label{B_dual}
L((\widehat{\ff})^{-1}) = \widehat{L({\ff}^{-1})} = [\beta_i - b_{ij}],
\end{equation}
where $L({\ff}^{-1})=[b_{ij}]$ and $\beta_i = \max \{b_{ij}\,|\, j=0, \ldots,n \}$.
\end{Proposition}
\proof
By Proposition~\ref{dual_mons_is_birational}, $\widehat{\ff}$ defines a Cremona map, so it remains to show
equality (\ref{B_dual}).
Consider the log-matrices $L(\ff)$ and $L(\ff^{-1})$.
By \cite[Theorem 2.2]{CremonaMexico}, there is a unique vector $\gamma\in \mathbb{N}^{n+1}$ (the inversion vector) such that
\begin{equation}\label{fund_relation}
L(\ff)\cdot L(\ff^{-1})=\Gamma+I_{n+1},
\end{equation}
where $\Gamma=[\underbrace{\gamma|\cdots |\gamma}_{n+1}]$.

Set $L(\ff)=(a_{ij})$ and $L(\ff^{-1})=(b_{ij})$ and let $ \alpha_i = \max\{a_{i0}, \ldots , a_{in}\}$
and $\beta_i = \max\{b_{i0}, \ldots , b_{i,n}\}$ denote the  $i$th coordinates
of the respective directrix vectors.

Define:
\begin{eqnarray}\label{inverse_gamma} \nonumber
\widehat{b_{ij}}: &=& \beta_i - b_{ij}\\ \nonumber
\widehat{\gamma_i}: &=& \gamma_i+ \alpha_i (\sum_{l=0}^{n} \beta_l - d') - \sum_{l=0}^{n} a_{il}\beta_l\\
&=& \gamma_i + \left(\alpha_i \deg(\widehat{\ff^{-1}})- (L(\ff)\cdot\underline{\beta}^t)_i \right).
\end{eqnarray}
where $d'=\deg(\ff^{-1})$ and $\underline{\beta}^t$ is the directrix vector of $\ff^{-1}$.
Note that $\widehat{\gamma_i}$ depends solely upon the $i$th rows of $L(\ff)$ and $L(\ff^{-1})$.

Letting $\widehat{c_{ij}}$ denote the $i,j$ entry of the product $\widehat{L(\ff)}\cdot (\widehat{b_{ij}})$, one has
\begin{eqnarray*}
\widehat{c_{ij}} &=& \sum_{l=0}^{n} \widehat{a_{il}} \widehat{b_{lj}}= \sum_{l=0}^{n} (\alpha_i-a_{il})(\beta_l - b_{lj})\\
&=& \sum_{l=0}^{n} \alpha_i\beta_l - \sum_{l=0}^{n} \alpha_i b_{lj} - \sum_{l=0}^{n} a_{il} \beta_l + \sum_{l=0}^{n} a_{il}b_{lj}.
\end{eqnarray*}
From this, using (\ref{fund_relation}), follows
\begin{equation*}
\widehat{c_{ij}} = \alpha_i(\sum_{l=0}^{n} \beta_l - d') - \sum_{l=0}^{n} a_{il}\beta_l + \gamma_i + \delta_{ij}
= \widehat{\gamma_i} + \delta_{ij}.
\end{equation*}
Therefore, setting $\widehat{\gamma} = (\widehat{\gamma_0}\ldots \widehat{\gamma_n})$ and
noting that $\widehat{L(\ff^{-1})} =[\widehat{b_{ij}}]$, one obtains
$$\widehat{L(\ff)} \widehat{L(\ff^{-1})} =
[\underbrace{\widehat{\gamma}^t|\cdots |\widehat{\gamma}^t}_{n+1}] + I_{n+1},$$
and hence, again by \cite[Theorem 2.2]{CremonaMexico}, $\widehat{L(\ff^{-1})}=L((\widehat{\ff})^{-1})$.
\qed

\begin{Remark}\rm Note that $\gamma$ and $\hat{\gamma}$ are additively related by a summand
that can be positive or negative.
It is worth observing that the monomial whose log-matrix is the vector $\gamma$ above
is exactly the so-called target inversion factor of the Cremona map defined by $\ff$.
This points to a potential meaningful relation between the respective inversion factors of $\ff$ and its Newton
complementary dual, at least in the monomial case.
Inversion factors play a role in the realm of symbolic powers (see \cite{Zaron}).
\end{Remark}

Unfortunately, Proposition~\ref{operations_commute} does not extend to arbitrary Cremona maps. The following is one of
the simplest examples where it fails.

\begin{Example}\rm
Consider the polar map in $\pp^2$ defined by the cubic form $f:=x(xz-y^2)\in k[x,y,z]$ (conic and a tangent line).
It is well known that this polynomial is homaloidal and the associated polar map is a de Jonqui\`eres
map of degree $2$, hence an involution (up to a change of variables both in the source and the target).

An easy calculation shows that the Newton complementary dual of the set of the partial derivatives $\{2xz-y^2,-2xy,x^2\}$
is $\{x^2z-2xy^2, xyz, y^2z\}$.
The claim is that the Cremona map defined by these cubic forms is not an involution
 (up to a change of variables both in the source and the target).
The inverse map to the the latter Cremona map is defined by the forms
$\{u^3-tuv, u^2v-tv^2, 2uv^2\}$ in the target variables $t,u,v$.
We cannot obviously argue by trying every change of variables in source and target. Thus, we need
an effective criterion. We use the one stated in \cite[Section 2.2]{AHA}, by which, in the particular case of
rational maps defined on the ambient $\pp^n$, there is a unique representative without fixed part.
Since the tuple $(x^2z-2xy^2, xyz, y^2z)$ has no fixed part, and neither does the tuple $(u^3-tuv, u^2v-tv^2, 2uv^2)$
even after change of variables and coordinates, the only way would be that they give the same tuple
up to renaming variables. This would imply, in particular, that the corresponding base ideals would
coincide. But this is not the case as the first has $3$ minimal primes, while the second has only $2$ minimal primes
(in particular, using a more geometric language, the maps have sets of proper and  infinitely
near points base points of different nature).
\end{Example}

\subsubsection{Plane de Jonqui\`eres maps}

The so-called de Jonqui\`eres maps are at the heart of classical plane Cremona map theory.
The next result enhances the role of the Newton dual for these maps.

We briefly recall a few preliminaries on these maps.
A plane de Jonqui\`eres map can be defined by $3$ forms of degree $d$ in $k[x,y,z]$ of the shape
$xq,yq,f$, where both $q$ and $f$ are $z$-monoids, with $f$ irreducible (see \cite[Section 4.1., especially Proposition 4.3]{HS}).
(A $z$-monoid is a polynomial of the shape $f=f_0+zf_1$, where $f_0,f_1$ are forms in $k[x,y]$.
Irreducibility of such an $f$ is equivalent to having $\gcd(f_0,f_1)=1$ -- see Section~\ref{sec:3}
for a more encompassing use of this notion.

The next result shows that de Jonqui\`eres maps are preserved under the Newton dual transform.

\begin{Proposition}\label{dual_of_dejonquieres}
The  Newton complementary dual of a plane de Jonqui\`eres map $\mathfrak{F}$ is a plane de Jonqui\`eres map.
More exactly, let $\mathfrak{F}$ be defined by nonzero forms $\ff=\{xq,yq,f\}\subset k[x,y,z]$ of degree $d$,
where $q=g+zh$, with $g,h\in k[x,y]$ of degrees $d-1, d-2$, respectively, and $f=f_0+zf_1$,
with $f_0,f_1\in k[x,y]$ of degrees $d, d-1$, respectively.
Then the linear system spanned by the Newton complementary dual $\widehat{\ff}$ defines a plane de Jonqui\`eres map
of degree
$$\hat{d}:= \max\{j_g+1, j_h+1, j_{f_0}, j_{f_1}\} - \min\{i_g, i_h+1, i_{f_0}, i_{f_1}+1\} +1,$$
where

\noindent$\bullet$  $0\leq i_g\leq d-1$ is the first index
such the term of $g$ of order $(d-1-i_g,i_g)$ does not vanish

\noindent $\bullet$ $0\leq j_g\leq d-1$  is the last index
such the term of $g$ of order $(d-1-j_g,j_g)$  does not vanish$\,${\rm ;}

\noindent $\bullet$
Similar definitions for $\{i_h, j_h\}$ {\rm (}respectively, $\{i_{f_0}, j_{f_0}\}$ and $\{i_{f_1} , j_{f_1}${\rm )}
regarding the form $h$ {\rm (}respectively, $f_0$ and $f_1${\rm )}.
\end{Proposition}
\proof
A straightforward inspection of the shape of the map gives the directrix vector
$$\underline{\alpha}=(\,\max\{d-i_g, d-1-i_h, d-i_{f_0}, d-i_{f_1}\} ,  \max\{j_g+1, j_h+1, j_{f_0}, j_{f_1}\},\, 1\,)$$
as defined earlier.
It remains to show that the resulting Cremona map is a de Jonqui\`eres map.

For this, we display the Newton matrices of the forms $\ff:=\{xq,yq,f\}$, respecting the lexicographic order:

\begin{equation}\label{xq}
N(xq)=\left(
\begin{array}{ccccc|ccccc}
d-i_g&>&\cdots &>& d-j_g&\;d-1-i_h &>& \cdots &>& d-1-j_h\\
i_g&<&\cdots &<&j_g&\;i_h &<&\cdots &<&j_h\\
0 && \cdots && 0 &\; 1 && \cdots && 1
\end{array} \right)
\end{equation}

\begin{equation}\label{yq}
N(yq)=\left(
\begin{array}{ccccc|ccccc}
d-1-i_g&>&\cdots &>& d-1-j_g&\;d-2-i_h &>& \cdots &>& d-2-j_h\\
i_g+1&<&\cdots &<&j_g+1&\;i_h+1 &<&\cdots &<&j_h+1\\
0 && \cdots && 0 &\; 1 && \cdots && 1
\end{array} \right)
\end{equation}

\begin{equation}\label{f}
N(f)=\left(
\begin{array}{ccccc|ccccc}
d-i_{f_0}&>&\cdots &>& d-j_{f_0}&\;d-1-i_{f_1}&>&\cdots &>& d-1-j_{f_1}\\
i_{f_0}&<&\cdots &<&j_{f_0}&\;i_{f_1}&<&\cdots &<&j_{f_1}\\
0 && \cdots && 0 &\; 1 && \cdots && 1
\end{array} \right)
\end{equation}

\smallskip

Here the symbol $>$  (respectively, $<$) is displayed to stress how the entries decrease
(respectively, increase).
Note that the sequences of entries represented by $\cdots$ are typically lacunary as many terms
have null coefficients. Actually, some blocks can collapse to one-column blocks.
As a slight control, note that $ d-j_g>0$ since, by definition, $d-1-i_g\neq 0$; etc.

Clearly, $\widehat{f}$ is again a $z$-monoid.
We next argue that two elements of $\widehat{\ff}$ are again of the form $xQ,yQ$
for some $z$-monoid $Q$.
For this we show that the first coordinate of $\alpha$ satisfies $\alpha_1 - N(yq)_{1i} \geq 1$ $\forall i$.
Namely, one has
$$\alpha_1 - N(yq)_{1i} \geq \alpha_1 - \max\{d-1-i_g, d-2-i_h\} \geq \max\{d-i_g, d-1-i_h\} - \max\{d-1-i_g, d-2-i_h\} =1.$$
One similarly shows that $\alpha_2 - N(yq)_{2i} \geq 1$ $\forall i$.
This implies that there are forms $Q_1,Q_2$ such that $\widehat{yq}=xQ_1$ and $\widehat{xq}=yQ_2$.
We claim that not only are $Q_1,Q_2$ $z$-monoids, but also $Q_1=Q_2$.
But this is immediate from looking at the respective Newton matrices
$$N(Q_1)=[\alpha - N(yq)] - \left(
\begin{array}{cccc}
1&1&\cdots &1 \\
0&0&\cdots &0\\
0 &0& \cdots &0 \end{array} \right) \quad {\rm and}\quad
N(Q_2)=[\alpha - N(xq)] - \left(
\begin{array}{cccc}
0&0&\cdots &0 \\
1&1&\cdots &1\\
0 &0& \cdots &0
 \end{array} \right).
 $$

To conclude, we know that since $\widehat{\ff}$ defines a Cremona map without fixed part, according to the criterion of
Enriques (\cite[Theorem 5.1.1]{alberich}) the $k$-vector space spanned by the elements of $\widehat{\ff}$
contains an irreducible form $F$.
Since the spanning elements are $z$-monoids, so must be $F$.
Clearly then this vector space is spanned by $\{xQ,yQ,F\}$, thus defining a de Jonqui\`eres map.
 \qed

\section{Cremona maps galore with prescribed invariants}
\label{sec:3}

We retain the terminology and notation of Subsection~\ref{sec:1.1}.
Thus, let $\mathfrak{F}:\pp^n\dasharrow \pp^n$ denote a Cremona map defined by a linear system spanned by
forms $f_0,\ldots, f_n$ of the same degree satisfying the usual canonical restrictions.
The common degree of these forms is the degree of $\mathfrak{F}$.
The ideal generated by these forms in the homogeneous coordinate ring $R=k[x_0,\ldots,x_n]$ of $\pp^n$
is the base ideal of $\mathfrak{F}$.

\begin{Definition}\rm We will say that a Cremona map has a certain property if its base ideal has
this property. Thus, a Cremona map whose base ideal is Cohen--Macaulay, of linear type, of codimension
$r$, etc. is a Cohen--Macaulay Cremona map,  a Cremona map of linear type,
of codimension $r$, and so forth.
\end{Definition}

The following elementary result will subsequently be used.

\begin{Lemma}\label{nzd}
Let $I\subset R$ be any ideal of a ring $R$, let $g\in R$ and let $X$ denote an indeterminate over $R$.
Then $Xg$ is a zerodivisor on $R{X}/IR[X]$ {\rm (}if and{\rm )} only if $g$ is a zerodivisor on $R/I$.
\end{Lemma}
\proof The ``if'' assertion is obvious.
To get the other direction, assume that $g$ is {\em not} a zerodivisor on $R/I$ and let
$h\in R[X]$ be such that $h\cdot Xg\in IR[X]$.
Then $hg\in IR[X]:(X)=IR[X]$ since $X$ is an indeterminate over $R$.
Therefore $h\in IR[X]:g= (I:g)R[X]$, i.e., every coefficient of $h$ belongs to $I:g$.
Since we are assuming that $g$ is not a zero divisors on $R/I$, it follows that $h\in IR[X]$, as was to be shown.
\qed

\begin{Remark}\rm The following souped-up version of Lemma~\ref{nzd} will be used below: {\em $Xg$ belongs to a minimal prime
of $R{X}/IR[X]$ of minimal codimension if and only if $g$ belongs to a minimal prime of $R/I$ of minimal
codimension}.
\end{Remark}

\subsection{Birational maps by recurrence}
\label{sec:3.1}

The  following recurrence proposition clarifies and extends the result proved in \cite[Proposition 3.2.2]{Barbara}.

\begin{Proposition}\label{Cremona_by_iteration}
Let $\mathfrak{F}:\pp^{n-1}\dasharrow \pp^{m-1}$ stand for a rational map defined by forms
$f_0,\ldots,f_{m-1}\in R=k[x_0,\ldots,x_{n-1}]$  of degree $d$.
Let
$$g:=g_{\ell}x_n^{d-\ell}+g_{\ell+1}x_n^{d-\ell-1}+\cdots +g_{d-1}x_n+g_d\in R[x_n]$$
be a form of degree $d$, where $g_j\in R$, for $\ell\leq j\leq d$, and $\deg_{x_n}(g)=d-\ell$, with $\ell\leq d-1$.
If $\mathfrak{F}$ is birational onto its image then the rational map $\mathfrak{G}:\pp^{n}\dasharrow \pp^{m}$
defined by the forms
$\{f_0, \ldots , f_{m-1},g\}\subset R[x_n]$ has topological {\rm (}i.e., field{\rm )} degree at most $d-\ell$.

In particular, if $\ell=d-1$ then $\mathfrak{G}$ is birational onto its image.

\smallskip

{\sc Supplement.} Suppose that $\deg_{x_n}(g)=1$ with $g_d\in I:=(f_0,\ldots,f_{m-1})$  {\rm (}possibly vanishing{\rm )}.
 If $g_{d-1}$ avoids all minimal primes of $R/I$ of minimal
codimension {\rm (}respectively, is a nonzerodivisor on $R/I$ and the latter is  Cohen--Macaulay{\rm )}
then the codimension of $R[x_n]/(I,g)$
is one larger than that of $R/I$ {\rm (}respectively, $R[x_n]/(I,g)$ is Cohen--Macaulay{\rm )}.
\end{Proposition}
\proof
By assumption, the field extension $k(f_1,\ldots,f_{m-1}) \subset  k((x_0,\ldots,x_{n-1})_d)$
is an equality, where $(x_0,\ldots,x_{n-1})_d$ denotes the set of monomials of degree $d$ in $k[x_0,\ldots,x_{n-1}]$.
On the other hand, it is easy to see that $k((x_0,\ldots,x_{n-1})_d)=k(x_0^d, x_1/x_0,\ldots,x_{n-1}/x_0)$.
Therefore, and by the same token, it suffices to prove that  $x_n/x_0$ satisfies an equation of
degree $d-\ell$ with coefficients in $k(f_0,\ldots,f_{m-1}, g)$.

Using the condensed monomial notation $\mathbf{x}^{\mathbf{a}}=x_0^{a_0}\cdots x_{n-1}^{a_{n-1}}$, let
$\sum_{|\mathbf{a}|=j} c_{\mathbf{a}}\mathbf{x}^{\mathbf{a}}$ be the expression of a form  $f\in R$ of degree $j$.
Then
\begin{eqnarray*} f &=& \sum_{|\mathbf{a}|=j} c_{\mathbf{a}} x_0^{\overbrace{a_0+\cdots+a_{n-1}}^{j}}
\underbrace{(x_1/x_0)^{a_1} \ldots
(x_{n-1}/x_0)^{a_{n-1}}}_{h_{\mathbf{a}} \in k(f_0,\cdots,f_{m-1})}\\
& = & x_0^{j} \sum_{|\mathbf{a}|=j} c_{\mathbf{a}} h_{\mathbf{a}}=x_0^{j} h_j,
\end{eqnarray*}
for suitable $h_j \in k(f_0,\cdots,f_{m-1})$,
where the equality  $k(f_0,\cdots,f_{m-1})=k(x_0^d, x_1/x_0,\ldots,x_{n-1}/x_0)$ given in the hypothesis
has been used.

 Therefore
 \begin{eqnarray*} g &=& x_0^{\ell}x_n^{d-\ell}h_{\ell} + x_0^{\ell+1}x_n^{d-\ell-1}h_{\ell+1} +\cdots + x_0^{d-1}x_nh_{d-1} + x_0^dh_d\\
 &=& x_0^{\ell}\left(x_n^{d-\ell}h_{\ell} + x_0x_n^{d-\ell-1}h_{\ell+1} +\cdots + x_0^{d-\ell-1}x_nh_{d-1} + x_0^{d-\ell}h_d\right)\\
 &=& x_0^d\left(h_{\ell}\left(\frac{x_n}{x_0}\right)^{d-\ell} + h_{\ell+1}\left(\frac{x_n}{x_0}\right)^{d-\ell-1}
 +\cdots + h_{d-1}\left(\frac{x_n}{x_0}\right) + h_d\right).
 \end{eqnarray*}
Now, $x_0^{-d} g\in k(x_0^d, x_1/x_0,\ldots,x_{n-1}/x_0,g)= k(f_0,\ldots,f_{m-1},g)$.
 Since $h_j\in k(f_0,\ldots,f_{m-1})$ for every $j$, we are through.

 \smallskip

The assertions in the supplement follow immediately from a slight modification of Lemma~\ref{nzd} and the remark after it,
using the assumption that $g_d\in I$.
\qed

\medskip

Forms of the shape $g=g_{d-1}x_n+g_d$, with $g_{d-1},g_d\in R$ are called {\em monoids}
 -- see \cite{Piene_et_al} for an overview of the concept.
It has been considered in Section~\ref{sec:1.4} in the case $n=2$ in connection to a de Jonqui\`eres map.
We will often refer to a monoid as an $x_{n}$-monoid in order to stress the privileged variable $x_{n}$.

\begin{Remark}\rm It is conceivable that the degree of $\mathcal{G}$ in Proposition~\ref{Cremona_by_iteration}
be exactly $d-\ell$ provided the coefficients $g_j$ are sufficiently general forms.
Drawing upon the characteristic-free methods developed in \cite{AHA}, where a certain Jacobian matrix is introduced,
it would amount to be able to express the rank of that matrix in a formula such as $\max\{0, n+1-(d-\ell)\}$.
Computation however can get pretty heavy.
Thus, even at the far end of the spectrum, when $d-\ell=1$ (the birational case) the degree of the inverse
map to $\mathcal{G}$ is not easily predictable in terms of the degree of the inverse to $\mathcal{F}$.
This tells us that the corresponding two graphs have far different structure -- alternatively,
the equations defining the Rees algebra of the ideal $(I,g)\subset R[x_n]$ are not totally predicted
from those defining the Rees algebra of the ideal $I\subset R$.
\end{Remark}

A slightly different type of construction yields Cremona maps which are Cohen--Macaulay with base ideal
of codimension $2$.
As above, we first state isolate a general result which is possibly available in
the literature in a disguised form.

\begin{Proposition}\label{CM_iteration}
Let $R$ be a Cohen--Macaulay Noetherian ring and let $I\subsetneq R$ be an ideal such that $R/I$
is perfect of codimension $2$.
Let $R\lar T$ be a flat $R$-algebra such that $IT\neq T$ and $\hht IT\geq \hht I$.
Then, for arbitrary elements $g\in IT$ and $f\in R$ such that the ideal $(f,g)T$ has codimension $2$,
the $T$-module $T/(If,g)T$ is  perfect of codimension $2$.
\end{Proposition}
\proof
Tensor a free $R$-resolution $0\rar F_1\stackrel{\phi}{\lar} F_0\lar R$ of $R/I$ with $T$ to get a similar
free $T$-resolution of $T/IT$. For simplicity, we write $\phi$ for $T\otimes_R\phi$ and $F_i$ for $T\otimes_R F_i$, $i=0,1$.

As one readily verifies, multiplication by $g$ induces an injective $T$-module homomorphism
$T/(IT:(g))f=T/IfT\injects T/IfT$ with image $(If,g)T/IfT$.
This homomorphism lifts to a map of complexes (free resolutions over $T$)
$$
\begin{array}{ccccccccccc}
0 & \rar & F_1 &\stackrel{\phi}{\lar} & F_0 & \stackrel{f\,\gg}{\lar} & T & \lar & T/IfT &\rar & 0\\[6pt]
 &&             \uparrow   &                 & \kern-12pt c(g)\uparrow   &    &\kern-8pt\cdot g\uparrow  && \cdot g\uparrow &&\\[6pt]
&& 0 &\rar & T     & \stackrel{f}{\lar} & T & {\lar} & T/IfT &\rar & 0
\end{array},
$$
where $\gg$ is the map induced by a set of generators of $I$  and $c(g):T\lar F_0$ is the content map that writes
$g$ as an element of $IT$, once a basis of $F_0$ is fixed.

Then the corresponding mapping cone is a $T$-free resolution (see \cite[Exercise A3.30]{E}):
\begin{equation}\label{mapping_cone}
0 \rar  F_1\oplus T \stackrel{\psi}{\lar}  F_0\oplus T  \stackrel{(f\,\gg, g)}{\lar}  T  \lar  T/(If,g)T \rar  0,
\end{equation}
where
$$\psi= \left(
\begin{array}{cc}
\phi & c(g)\\
\mathbf{0} & -f
\end{array}
\right).
$$
To conclude, the ideal $(If, g)\subset T$ has codimension $2$ since, by assumption, $(f,g)\subset T$ has
codimension $2$ and $(I,g)T=IT$ has codimension $\geq 2$.
\qed

\bigskip

The proposition has the following application to birational theory.

\begin{Theorem}\label{CM_birational_cod2}
Let $f_0,\ldots,f_{m-1}\in R=k[x_0,\ldots,x_{n-1}]$ be forms of degree $d$ satisfying the canonical
restrictions and let $\mathfrak{F}:\pp^{n-1}\dasharrow \pp^{m-1}$ denote the corresponding rational map.
Set $I:=(f_0,\ldots,f_{m-1})\subset R$.
Let $f\in R$ and $g\in IR[x_n]$, be nonzero forms such that $\deg(g)=\deg(f)+d$ and $\gcd(f,g)=1$
in $R[x_n]$, where  $x_n$ is a new variable.
Suppose that:
\begin{itemize}
\item The base ideal $I$ is a codimension $2$ perfect ideal.
\item $g$ is an $x_n$-monoid.
\item $\mathfrak{F}$ is birational onto its image {\rm (}hence, $m\geq n${\rm )}.
\end{itemize}
\noindent Then the rational map $\mathfrak{G}:\pp^{n}\dasharrow \pp^{m}$ defined by the forms
$\{f_0f,\ldots , f_{m-1}f,g\}\subset R[x_n]$
is a birational map {\rm (}onto its image{\rm )} whose base ideal is perfect of codimension $2$.
In particular, if $\mathfrak{F}$ is a perfect codimension $2$ Cremona map of $\pp^{n-1}$
then $\mathfrak{G}$ is a perfect codimension $2$ Cremona map of $\pp^{n}$.
\end{Theorem}
\proof
Note that the assertion can be broken up into two separate assertions.
The assertion about the perfectness of the base ideal $(If,g)$ follows immediately from Proposition~\ref{CM_iteration}
by taking $T:=R[x_n]$ and applying it to the ideal $I=(f_0,\ldots,f_{m-1})$.

To prove the birationality claim one argues that $f_0f,\ldots , f_{m-1}f$ also defines a birational map,
hence one can apply Proposition~\ref{Cremona_by_iteration} to this map ({\em not} to the rational map
defined by $f_0,\ldots,f_{m-1}$, and observe we do not use the supplement there).
\qed

\begin{Remark}\rm The content of the above theorem is similar to the situation of a {\em de Jonqui\`eres parametrization}
as in \cite{stellar} -- actually, to a special case called ``inclusion case''. However, the present setup is slightly
different as one is feeding in a new variable as well (a flat extension).
This has the effect that the proof in the present situation is a bit less automatic.
Moreover, as easy examples show, the hypothesis of $g$ being a monoid is in general essential both in
Theorem~\ref{CM_birational_cod2} and  Proposition~\ref{Cremona_by_iteration}.
\end{Remark}

The following is the main result of the section.

\begin{Theorem}\label{arbitrary_cod}
Let $n\geq 2,\, d\geq 2$ be  integers. Then, for every $r$ in the integer interval
$[2,n]$,  there exist Cremona transformations of $\pp^n$ of codimension $r$ and degree $d$,
and there exist Cohen--Macaulay Cremona transformations of $\pp^n$ of codimension $r$ and
degree $d\geq n+1-r$.
\end{Theorem}
\demo Induct on $n$.
If $n=2$ there is nothing to prove since any Cremona transformation of $\pp^2$ has  codimension $2$
and, for any $d\geq 1$ the forms $x_0^d, x_0^{d-1}x_1, x_1^{d-1}x_2$ define a plane Cohen--Macaulay Cremona map of degree $d$.

So, suppose that $n\geq 3$. By the inductive hypothesis, for any $r$ in the integer range $[2,n-1]$,
 $\pp^{n-1}$ admits  Cremona maps of codimension $r$ and degree $d$, and also
Cohen--Macaulay ones of of codimension $r$ and degree $d\geq n-1+1-r = n+1-(r+1)$.
By Proposition~\ref{Cremona_by_iteration}, for any integer $s$ in the range $[2,n]$, we obtain
Cremona maps of codimension $s$ and degree $d$, and also
Cohen--Macaulay Cremona transformations of $\pp^n$ of any codimension in the range $[3,n]$ and degree $d\geq n+1-s$.
To show the existence of the latter also in codimension $2$ we apply Theorem~\ref{CM_birational_cod2}
by initiating with any Cohen--Macaulay plane Cremona (see the remark below and also
Theorem~\ref{sub_Hankel} in the next section which gives a structured example of a perfect codimension $2$ Cremona map of
$\pp^{n-1}$ of any degree $\geq n-2+q$ with $q\geq 1$.)
\qed

\begin{Remark}\rm
For $n=2$ and arbitrary $d\geq 1$, the class of de Jonqui\`eres maps gives infinitely many distinct Cohen--Macaulay examples
(see, e.g., \cite[Proposition 4.2 and Corollary 4.5]{HS}), of which the monomial Cremona map in the proof is an instance.
Thus,  the proof of the theorem also shows that, for any $n\geq 2$, there are infinitely many such
Cohen--Macaulay Cremona maps for any choices of
$d\geq 1$ and of a prescribed codimension in the integer range $[2,n]$.
\end{Remark}

The lower bound  $d\geq n+1-r$, with $r$ the codimension of a Cohen--Macaulay Cremona map,
is quite natural. Actually, it always holds for $r=2$ since the lowest (standard) degree of a syzygy of the base ideal
is $1$, while it adds no condition for $r=n$.
The natural question left is:

\begin{Question}\rm Is the lower bound $d\geq n+1-r$, in the range $r\in [3,n-1]$, for Cohen--Macaulay Cremona maps,
best possible for all $n\geq 3$?
\end{Question}

\subsection{A determinantal model}
\label{sec:3.2}

The following degeneration of the ordinary generic
Hankel matrix has been introduced in {\rm \cite[Section 4.1.1]{CRS}}, where it has been called a {\em sub-Hankel} matrix.

\begin{equation}\label{subH}
\left(
\begin{matrix}
x_0&x_1&x_2&...&x_{n-2}&x_{n-1}\\
x_1&x_2&x_3&...&x_{n-1}&x_{n}\\
x_2&x_3&x_4&...&x_{n}&0\\
.&.&.&...&.&.\\
.&.&.&...&.&.\\
.&.&.&...&.&.\\
x_{n-2}&x_{n-1}&x_{n}&...&0&0\\
x_{n-1}&x_{n}&0&...&0&0\\
\end{matrix}
\right)
\end{equation}

It has been shown in \cite{CRS} that the determinant $F\in k[x_0,\ldots,x_n]$ of the sub-Hankel matrix
is a homaloidal polynomial (with degree equal to the dimension of the ambient projective space).
As a bonus one also obtained Cremona maps with Cohen--Macaulay codimension $2$ base ideal.
However, the degree of any of these maps is bounded above by the number of variables.
Here we wish to use part of the coordinates of the corresponding polar map -- i.e., a subset of the partial derivatives of
$F$ -- in order to obtain a determinantal like model for the Cremona maps in the previous subsection which also
have Cohen--Macaulay codimension $2$ base ideal and, moreover, of degree arbitrarily
larger than the number of variables.

\begin{Theorem}\label{sub_Hankel}
Consider the  matrix {\rm (\ref{subH})}, with $n\geq 3$ and determinant $F\in k[x_0,\ldots,x_n]$.
Let  the partial derivatives of $F$ be denoted $F_{x_0},\ldots , F_{x_{n-1}}, F_{x_n}$.
Then $x_n$ divides $F_{x_i}$, for $0\leq i\leq n-2$ and for
any integer $q\geq 1$, the forms
\begin{equation*}
\frac{F_{x_0}}{x_n}\,x_n^{q},\,\ldots\, , \frac{F_{x_{n-2}}}{x_n}\,x_n^{q},\;
F_{x_{n-1}}x_n^{q-1} + \frac{F_{x_{n-2}}}{x_n}\,x_{n-1}^{q}
\end{equation*}
belong to $k[x_1,\ldots,x_n]$ and define a codimension $2$ Cohen--Macaulay Cremona map of $\,\pp^{n-1}$.
\end{Theorem}
\demo
First, one has $\gcd(F_{x_0},\ldots,F_{x_{n-2}})=x_n$ and $x_n$ does not divide $F_{x_{n-2}}/x_n$ (second claim
of \cite[Lemma 4.2 (ii)]{CRS}).

Moreover, The first $n-2$ partial derivatives of $F$ divided by their gcd $x_n$ define a Cremona map
with perfect codimension $2$ base ideal (\cite[Theorem 4.4 (i)]{CRS}).

Next, $F_{x_0},\ldots,F_{x_{n-2}}\in k[x_2,\ldots,x_n]$ (\cite[Lemma 4.2 (ii)]{CRS}).
Also, by direct inspection of the derivatives one sees that $F_{x_{n-i}}$ is an $x_i$-monoid.
In particular, $F_{x_{n-1}}$ is an $x_1$-monoid, hence $g:=F_{x_{n-1}}x_n^{q-1} + \frac{F_{x_{n-2}}}{x_n}\,x_{n-1}^{q}$
is also an $x_1$-monoid since $x_1$ does not appear in $F_{x_{n-2}}$.
It is clear that $g$ is not divisible by $x_n$.

Finally, $g$ belongs to the ideal generated by the first $n-2$ partial derivatives of $F$ divided
 by their gcd $x_n$ as a case of the relation (4.4) in \cite[Lemma 4.2 (iii)]{CRS}.

Assembling the information, the result is now an application of Theorem~\ref{CM_birational_cod2},
with $f:=x_n^q$ and $f_i:=F_{x_i}/x_n,\, i=0,\ldots, n-2$ and $x_1$ as the additional variable
(in the role of $x_n$ in that theorem).
\qed

\begin{Remark}\rm Of course, one could trade the last coordinate $g$
 in the theorem by any $x_1$-monoid of the right degree not divisible by $x_n$ and belonging to
 the ideal generated by the first $n-2$ partial derivatives of $f$ divided
 by their gcd $x_n$.
 However, we wished to keep the map as close as possible to the original polar map.
 Moreover, there is some  computational evidence for conjecturing that its inverse map
has the same shape, though of a different degree.
\end{Remark}

\medskip

{\small
\noindent {\bf Authors' addresses:}

\medskip

\noindent {\sc Barbara Costa}\\
 Departamento de Matem\'atica, Universidade Federal Rural de Pernambuco \\
             52171-900, Recife, PE, Brazil\\
             {\em e-mail}: binhamat@gmail.com           
 
 \medskip
         
\noindent {\sc Aron Simis}\\  Departamento de Matem\'atica, CCEN, Universidade Federal
de Pernambuco\\
 50740-560 Recife, PE, Brazil.\\
{\em e-mail}:  aron@dmat.ufpe.br

}

\end{document}